\input amstex
\magnification=\magstep1 
\baselineskip=13pt
\documentstyle{amsppt}
\vsize=8.7truein \CenteredTagsOnSplits \NoRunningHeads
\def\UU{\Cal U}
\def\per{\operatorname{per}}
\def\haf{\operatorname{haf}}
\def\PER{\operatorname{PER}}
 \topmatter
 
\title Computing the permanent of (some) complex matrices   \endtitle 
\author Alexander Barvinok \endauthor
\address Department of Mathematics, University of Michigan, Ann Arbor,
MI 48109-1043, USA \endaddress
\email barvinok$\@$umich.edu \endemail
\date June 2014 \enddate
\thanks  This research was partially supported by NSF Grant DMS 0856640.
\endthanks 
\abstract We present a deterministic algorithm, which, for any given $0< \epsilon < 1$ and an $n \times n$ real or complex matrix 
$A=\left(a_{ij}\right)$ such that $\left| a_{ij}-1 \right| \leq 0.19$ for all $i, j$ computes the permanent of $A$ within relative error $\epsilon$ in $n^{O\left(\ln n -\ln \epsilon\right)}$ time. The method can be extended to computing hafnians and multidimensional permanents. \endabstract

\endtopmatter

\document

\head 1. Introduction and main results \endhead

The {\it permanent} of an $n \times n$ matrix $A=\left(a_{ij}\right)$ is defined as 
$$\per A =\sum_{\sigma \in S_n} \prod_{i=1}^n a_{i \sigma(i)},$$
where $S_n$ is the symmetric group of permutations of the set $\{1, \ldots, n\}$. The problem of efficient computation of the permanent has attracted a lot of attention. It is $\#P$-hard already for 0-1 matrices  \cite{Va79}, but a fully polynomial randomized approximation scheme, based on the Markov Chain Monte Carlo approach, is constructed for all non-negative matrices \cite{J+04}. A deterministic polynomial time algorithm based on matrix scaling for computing the permanent of non-negative matrices within a factor of $e^n$ is constructed in \cite{L+00} and the bound was recently improved to $2^n$ in
\cite{GS13}. An approach based on the idea of ``correlation decay" from statistical physics results in a deterministic polynomial time algorithm approximating $\per A$ within a factor of $(1+\epsilon)^n$ for any $\epsilon >0$, fixed in advance, if $A$ is the adjacency matrix of a constant degree expander \cite{GK10}.

There is also interest in computing permanents of {\it complex} matrices \cite{AA13}. The well-known Ryser's algorithm (see, for example, Chapter 7 of \cite{Mi78}) computes the permanent of a matrix $A$ over any field in $O\left(n 2^n\right)$ time. A randomized approximation algorithm of \cite{F\"u00} computes the permanent of a complex matrix within a (properly defined) relative error $\epsilon$ in $O\left(3^{n/2}\epsilon^{-2}\right)$ time. The randomized algorithm of \cite{Gu05}, see also \cite{AA13} for an exposition, computes the permanent of a complex matrix $A$ in polynomial in $n$ and $1/\epsilon$ time within an additive error of $\epsilon \|A\|^n$, where $\|A\|$ is the operator norm of $A$.

In this paper, we present a new approach to computing permanents of real or complex matrices $A$ and show that if $\left| a_{ij}-1\right| \leq \gamma$ for some absolute constant
$\gamma > 0$ (we can choose $\gamma=0.19$) and all $i$ and $j$, then, for any $\epsilon > 0$ the value of $\per A$ can be computed within relative error $\epsilon$ in $n^{O\left(\ln n- \ln \epsilon\right)}$ time (we say that $\alpha \in {\Bbb C}$ approximates 
$\per A$ within relative error $0< \epsilon < 1$ if $\per A =\alpha(1+ \rho)$ where
$|\rho|  < \epsilon$). We also discuss how the method can be extended to computing hafnians of symmetric matrices and multidimensional permanents of tensors.

\subhead (1.1) The idea of the algorithm \endsubhead
Let $J$ denote the $n \times n$ matrix filled with 1s. Given an $n \times n$ complex matrix $A$, we consider (a branch of) the univariate function 
$$f(z)=\ln \per \bigl(J + z(A-J)\bigr). \tag1.1.1$$
Clearly,
$$f(0)=\ln \per J =\ln n! \quad \text{and} \quad f(1)=\ln \per A.$$
Hence our goal is to approximate $f(1)$ and we do it by using the Taylor polynomial expansion of $f$ at $z=0$:
$$f(1) \approx f(0) +\sum_{k=1}^m {1 \over k!} {d^k \over dz^k} f(z) \Big|_{z=0}. \tag1.1.2$$
It turns out that the right hand side of (1.1.2) can be computed in $n^{O(m)}$ time. We present the algorithm in Section 2. The quality of the approximation (1.1.2) depends on the location of complex zeros of the permanent.

\proclaim{(1.2) Lemma} Suppose that there exists a real $\beta>1$ such that 
$$\per \bigl(J +z(A-J)\bigr) \ne 0 \quad \text{for all} \quad z \in {\Bbb C} \quad \text{satisifying} \quad |z| \leq \beta.$$
Then for all $z \in {\Bbb C}$ with $|z| \leq 1$ the value of 
$$f(z) =\ln \per \bigl(J+ z(A-J)\bigr)$$ is well-defined by the choice of the branch of the logarithm for which $f(0)$ is a real number, and the right hand side of (1.1.2) approximates $f(1)$ within an additive error of 
$${n\over (m+1)\beta^m(\beta-1)}.$$
\endproclaim
In particular, for a fixed $\beta >1$, to ensure an additive error of $0< \epsilon < 1$, we can choose 
$m=O\left(\ln n -\ln \epsilon\right)$, which results in the algorithm for approximating $\per A$ within relative error $\epsilon$ in
$n^{O\left(\ln n -\ln \epsilon\right)}$ time. We prove Lemma 1.2 in Section 2.

Thus we have to identify a class of matrices $A$ for which the number $\beta >1$ of Lemma 1.2 exists.
We prove the following result.

\proclaim{(1.3) Theorem} There is an absolute constant $\delta >0$ (we can choose $\delta=0.195$) such that if 
$Z=\left(z_{ij}\right)$ is a complex $n \times n$ matrix satisfying 
$$\left|z_{ij}-1 \right| \leq \delta \quad \text{for all} \quad i,j$$
then
$$\per Z \ne 0.$$
\endproclaim

We prove Theorem 1.3 in Section 3.

For any matrix $A=\left(a_{ij}\right)$ satisfying 
$$\left| a_{ij} -1 \right| \leq 0.19 \quad \text{for all} \quad i,j,$$
we can choose $\beta=195/190$ in Lemma 1.2 and thus obtain an approximation algorithm for computing $\per A$.

The sharp value of the constant $\delta$ in Theorem 1.3 is not known to the author. A simple example of a $2 \times 2$ matrix 
$$A=\left(\matrix {1+i \over 2} & {1-i \over 2}  \\ {1-i \over 2} & {1+i \over 2}  \endmatrix \right)$$
for which $\per A=0$ shows that in Theorem 1.3 we must have 
$$\delta \ < \ {\sqrt{2} \over 2} \approx 0.71.$$ 
What is also not clear is whether the constant $\delta$ can {\it improve} as the size of the matrix grows.

\subhead (1.4) Question \endsubhead Is it true that for any $0 < \epsilon < 1$ there is a positive integer $N(\epsilon)$ such that 
if $Z=\left(z_{ij}\right)$ is a complex $n \times n$ matrix with $n > N(\epsilon)$ and 
$$\left| z_{ij}-1\right| \ \leq \ 1-\epsilon \quad \text{for all} \quad i,j$$
then $\per Z \ne 0$?

We note that for any $0 < \epsilon < 1$, fixed in advance, a deterministic polynomial time algorithm based on scaling approximates the permanent of a given $n \times n$ real matrix $A=\left(a_{ij}\right)$ satisfying
$$\epsilon \ \leq \ a_{ij} \ \leq \ 1 \quad \text{for all} \quad i, j$$
within a multiplicative factor of $n^{\kappa(\epsilon)}$ for some $\kappa(\epsilon)>0$ \cite{BS11}.

\subhead (1.5) Ramifications \endsubhead In Section 4, we discuss how our approach can be used for computing hafnians of symmetric matrices and multidimensional permanents of tensors. The same approach can be used for computing partition functions associated with cliques in graphs \cite{Ba14} and graph homomorphisms \cite{BS14}.
In each case, the main problem is to come up with a version of Theorem 1.3 bounding the complex roots of the partition function away from the vector of all 1s. Isolating zeros of complex extensions of real partition functions is a problem studied in statistical physics and also in connection to combinatorics, see, for example, \cite{SS05}.  

\head 2. The algorithm \endhead

\subhead (2.1) The algorithm for approximating the permanent \endsubhead Given an $n \times n$ complex matrix
$A=\left(a_{ij}\right)$, we present an algorithm which computes the right hand side of the approximation (1.1.2) for the function $f(z)$ defined by (1.1.1).

Let 
$$g(z)=\per \bigl(J+z(A-J)\bigr), \tag2.1.1$$
so $f(z)=\ln g(z)$.
Hence
$$f'(z)={g'(z) \over g(z)} \quad \text{and} \quad g'(z)=g(z) f'(z).$$
Therefore, for $k \geq 1$ we have 
$${d^k \over dz^k}g(z)\Big|_{z=0} =\sum_{j=0}^{k-1} {k-1 \choose j} \left( {d^j \over dz^j} g(z) \Big|_{z=0} \right) 
\left({d^{k-j} \over dz^{k-j}} f (z) \Big|_{z=0}\right) \tag2.1.2$$
(we agree that the 0-th derivative of $g$ is $g$).

We note that $g(0)=n!$. If we compute the values of 
$${d^k \over dz^k} g(z)\Big|_{z=0} \quad \text{for} \quad k=1, \ldots, m, \tag2.1.3$$
then the formulas (2.1.2) for $k=1, \ldots, m$ provide a non-degenerate triangular system of linear equations that allows us to 
compute 
$${d^k \over dz^k} f(z)\Big|_{z=0} \quad \text{for} \quad k=1, \ldots, m.$$
Hence our goal is to compute the values (2.1.3).

We have 
$$\split {d^k \over dz^k} g(z)\Big|_{z=0}= &{d^k \over dz^k} \sum_{\sigma \in S_n}
\prod_{i=1}^n \bigl(1 + z\left(a_{i\sigma(i)}-1\right)\bigr) \Big|_{z=0}\\=
&\sum_{\sigma \in S_n} \sum_{1 \leq i_1, \ldots, i_k \leq n} \left(a_{i_1 \sigma(i_1)}-1\right) \cdots \left(a_{i_k \sigma(i_k)} -1 \right) \\
=&(n-k)!
 \sum \Sb 1 \leq i_1, \ldots, i_k \leq n \\ 1\leq j_1, \ldots, j_k \leq n \endSb \left(a_{i_1 j_1}-1\right) \cdots \left(a_{i_k j_k} -1 \right),
 \endsplit$$
where the last sum is over all pairs of ordered $k$-subsets $\left(i_1, \ldots, i_k \right)$ and $\left(j_1, \ldots, j_k \right)$
of the set $\{1, \ldots, n\}$. Since the last sum contains $\bigl(n!/(n-k)!\bigr)^2=n^{O(k)}$ terms, the complexity of the algorithm is indeed $n^{O(m)}$.

\subhead (2.2) Proof of Lemma 1.2 \endsubhead The function $g(z)$ defined by (2.1.1) is a polynomial in $z$ of degree 
$d \leq n$ with $g(0)=n! \ne 0$, so we factor
$$g(z)= g(0) \prod_{i=1}^d \left(1-{z \over \alpha_i}\right)
,$$
 $\alpha_1, \ldots, \alpha_d$ are the roots of $g(z)$. By the condition of Lemma 1.2, 
we have 
$$\left|\alpha_i\right| \ \geq \ \beta >1 \quad \text{for} \quad i=1, \ldots, d.$$
Therefore,
$$f(z)=\ln g(z)= \ln g(0) + \sum_{i=1}^d \ln \left(1-{z \over \alpha_i}\right) \quad \text{for} \quad |z| \leq 1, \tag2.2.1$$
where we choose the branch of $\ln g(z)$ that is real at $z=0$. Using the standard Taylor expansion, we obtain 
$$\ln \left(1-{1 \over \alpha_i}\right) =-\sum_{k=1}^m {1 \over k} \left({1 \over \alpha_i}\right)^k + \zeta_m,$$
where
$$\left|\zeta_m\right|= \left| \sum_{k=m+1}^{+\infty} {1 \over k} \left( {1 \over \alpha_i} \right)^k \right| \ \leq \ 
{1 \over (m+1) \beta^m (\beta-1)}.$$
Therefore, from (2.2.1) we obtain
$$f(1)=f(0)  + \sum_{k=1}^m \left(-{1 \over k} \sum_{i=1}^d \left({1 \over \alpha_i}\right)^k\right)  +\eta_m,$$
where 
$$\left| \eta_m \right| \ \leq \ {n \over (m+1) \beta^m (\beta-1)}.$$
It remains to notice that 
$$-{1 \over k} \sum_{i=1}^d \left({1 \over \alpha_i}\right)^k ={1 \over k!} {d^k \over dz^k} f(z) \Big|_{z=0}.$$
{\hfill \hfill \hfill} \qed

\head 3. Proof of Theorem 1.3 \endhead

Let us denote  by $\UU^{n \times n}(\delta) \subset {\Bbb C}^{n \times n}$ the closed polydisc 
$$\UU^{n \times n}(\delta) =\Bigl\{ Z=\left(z_{ij} \right): \quad \left| z_{ij}-1 \right| \ \leq \ \delta \quad \text{for all} \quad i, j \Bigr\}.$$
Thus Theorem 1.3 asserts that $\per Z \ne 0$ for $Z \in \UU^{n \times n}(\delta)$ and $\delta=0.195$.

First, we establish a simple geometric lemma.

\proclaim{(3.1) Lemma} Let $u_1, \ldots, u_n \in {\Bbb R}^d$ be non-zero vectors such that for some $0 \leq \alpha < \pi/2$
the angle between any two vectors $u_i$ and $u_j$ does not exceed $\alpha$. Let 
$u=u_1 + \ldots + u_n$. Then 
$$\|u\| \ \geq \ \sqrt{\cos \alpha}  \sum_{i=1}^n \|u_i\|.$$
\endproclaim
\demo{Proof} We have 
$$\|u\|^2 =\sum_{1 \leq i, j \leq n} \langle u_i, u_j \rangle \ \geq \ \sum_{1 \leq i, j \leq n} \|u_i\| \|u_j\| \cos \alpha =
\left(\cos \alpha \right) \left( \sum_{i=1}^n \|u_i \| \right)^2,$$
and the proof follows.
{\hfill \hfill \hfill} \qed
\enddemo

We prove Theorem 1.3 by induction on $n$, using Lemma 3.1 and the following two lemmas.

\proclaim{(3.2) Lemma} For an $n \times n$ matrix $Z=\left(z_{ij}\right)$ and $j=1, \ldots, n$, let $Z_j$ be the $(n-1) \times (n-1)$ matrix obtained from $Z$ by crossing out the first row and the $j$-th column of $Z$. 

Suppose for some $\delta > 0$ and for some $0 < \tau  < 1$, for any $Z \in \UU^{n \times n}(\delta)$ we have $\per Z \ne 0$ and
$$\left| \per Z \right| \geq \ \tau \sum_{j=1}^n \left| z_{1j}\right| \left| \per Z_j \right|.$$

Let $A, B \subset \UU^{n \times n}(\delta)$ be any two $n \times n$ matrices that differ in one column (or in one row) only.
Then the angle between two complex numbers $\per A$ and $\per B$, interpreted as vectors in ${\Bbb R}^2={\Bbb C}$ does not exceed 
$$\theta ={2 \delta \over (1-\delta) \tau}.$$
\endproclaim
\demo{Proof} Since $\per Z \ne 0$ for all $Z \in \UU^{n \times n}(\delta)$, we may consider a branch of $\ln \per Z$ defined for $Z \in \UU^{n \times n}(\delta)$.

Using the expansion 
$$\per Z =\sum_{j=1}^n z_{1j} \per Z_j, \tag3.2.1$$
we conclude that 
$${\partial \over \partial z_{1j}} \ln \per Z = {\per Z_j \over \per Z} \quad \text{for} \quad j=1, \ldots, n.$$
Therefore, since $\left|z_{ij}\right| \geq 1-\delta$ for $j=1, \ldots, n$, we conclude that for any $Z \in \UU^{n \times n}(\delta)$, we have
$$\sum_{j=1}^n \left| {\partial  \over \partial z_{1j}} \ln \per Z \right| \ \leq \ {1 \over (1-\delta) \tau}. \tag3.2.2$$
Since the permanent is invariant under permutations of rows, permutations of columns and taking the transpose of the matrix,
without loss of generality we may assume that the matrix $B \in \UU^{n \times n}(\delta)$ is obtained from $A \in \UU^{n \times n}(\delta)$ by replacing the entries $a_{1j}$ by numbers $b_{1j}$ such that 
$$\left| b_{1j}-1 \right| \leq \delta \quad \text{for} \quad j=1, \ldots, n.$$
Then 
$$\left| \ln \per A - \ln \per B \right| \ \leq \ \left( \sup_{Z \in \UU^{n \times n}(\delta)}
\sum_{j=1}^n \left| {\partial  \over \partial z_{1j}} \ln \per Z \right| \right) \left( \max_{j=1, \ldots, n} \left| a_{1j} -b_{1j} \right| \right).$$
Since 
$$\left| b_{1j} - a_{1j} \right| \ \leq \ 2 \delta \quad \text{for all} \quad j=1, \ldots, n,$$
the proof follows from (3.2.2).
{\hfill \hfill \hfill} \qed
\enddemo

\proclaim{(3.3) Lemma} Suppose that for some 
$$0\ \leq\ \theta \ < \ {\pi \over 2} -2 \arcsin \delta$$
and for any two matrices $A, B \in \UU^{n \times n}(\delta)$ which differ in one row (or in one column), the angle between two complex numbers $\per A$ and $\per B$, interpreted as vectors in ${\Bbb R}^2={\Bbb C}$ does not exceed $\theta$. Then for any matrix $Z \in \UU^{(n+1)\times (n+1)}(\delta)$, we have 
$$\left| \per Z \right| \ \geq \ \tau \sum_{j=1}^{n+1} \left|z_{1j}\right| \left| \per Z_j\right|$$
with 
$$\tau=\sqrt{\cos\left(\theta + 2 \arcsin \delta\right)},$$
where $Z_j$ is the $n\times n$ matrix obtained from $Z$ by crossing out the first row and the $j$-th column.
\endproclaim
\demo{Proof} We use the first row expansion (3.2.1) and observe that any two matrices $Z_j$ and $Z_k$, can be obtained from one from another by a replacing one column and a permutation of columns.
Therefore, the angle between any two complex numbers $\per Z_j$ and $\per Z_k$ does not exceed $\theta$. Since 
$$-\arcsin \delta \ \leq \ \arg z_{1j} \ \leq \ \arcsin \delta \quad \text{for} \quad j=1, \ldots, n,$$
the angle between any two numbers $z_{1j} \per Z_j$ and $z_{1k} \per Z_k$ does not exceed 
$\theta + 2 \arcsin \delta$. 
The proof follows by Lemma 3.1.
{\hfill \hfill \hfill} \qed
\enddemo

\subhead (3.4) Proof of Theorem 1.3 \endsubhead One can see that for a sufficiently small $\delta>0$, the equation 
$$\theta = {2 \delta \over (1-\delta)\sqrt{\cos(\theta +2 \arcsin \delta)}} \tag3.4.1$$
has a solution $0< \theta < \pi/2$. Numerical computations show that we can choose $\delta=0.195$ and 
$$\theta \approx 0.7611025127.$$
Let
$$\tau=\sqrt{\cos(\theta+2 \arcsin \delta)} \approx 0.6365398112.$$

We proceed by induction on $n$. More precisely, 
we prove the following three statements (3.4.2)--(3.4.4) by induction on $n$:
\bigskip
(3.4.2) For every $Z \in \UU^{n \times n}(\delta)$, we have $\per Z \ne 0$;
\medskip
(3.4.3) Suppose $A, B \in \UU^{n \times n}(\delta)$ are two matrices which differ by one row (or one column). Then the angle between two complex numbers $\per A$ and $\per B$, interpreted as vectors in ${\Bbb R}^2 ={\Bbb C}$, does not exceed $\theta$;
\medskip
(3.4.4) For a matrix $Z \in \UU^{n \times n}(\delta)$, $Z=\left(z_{ij}\right)$, let $Z_j$ be the $(n-1) \times (n-1)$ matrix obtained by crossing out the first row and the $j$-th column. Then 
$$\left|\per Z\right| \ \geq \ \tau \sum_{j=1}^n \left|z_{1j}\right| \left| \per Z_j \right|.$$
\bigskip
For $n=1$ the statement (3.4.2) is obviously true. Moreover, the angle between any two numbers $a, b \in \UU^{1 \times 1}(\delta)$ does not exceed 
$$2 \arcsin \delta \approx 0.3925149004 < \theta,$$
 so (3.4.3) holds as well. 
The statement (3.4.4) is vacuous.

Lemma 3.3 implies that if the statement (3.4.3) holds for $n \times n$ matrices then the statement (3.4.4) holds for 
$(n+1) \times (n+1)$ matrices. 

The statement (3.4.4) for $(n+1)\times (n+1)$ matrices together with the statement (3.4.2) for $n \times n$ matrices implies the statement (3.4.2) for $(n+1)\times (n+1)$ matrices.

Finally, Lemma 3.2 implies that if the statement (3.4.4) holds for $(n+1) \times (n+1)$ matrices then the statement (3.4.3) holds for $(n+1) \times (n+1)$ matrices.

This concludes the proof of (3.4.2)--(3.4.4) for all positive integer $n$.
{\hfill \hfill \hfill} \qed

\head 4. Ramifications \endhead

A similar approach can be applied to computing other quantities of interest.

\subhead (4.1) Hafnians \endsubhead Let $A=\left(a_{ij}\right)$ be a $2n \times 2n$ real or complex matrix. The quantity
$$\haf A=\sum_{\{i_1, j_1\}, \ldots, \{i_n, j_n\}} a_{i_1 j_1} \cdots a_{i_n j_n},$$
where sum is taken over all $(2n)!/n! 2^n$ unordered partitions of the set $\{1, \ldots, 2n\}$ into $n$ pairwise disjoint unordered pairs 
$\{i_1, j_1\}, \ldots, \{i_n, j_n\}$, is called the {\it hafnian} of $A$, see for example, Section 8.2 of \cite{Mi78}.
For any $n \times n$ matrix $A$ we have 
$$\haf\left(\matrix 0 & A \\ A^T & 0 \endmatrix \right)=\per A$$
and hence computing the permanent of an $n \times n$ matrix reduces to computing the hafnian of a symmetric $2n \times 2n$ matrix. The computational complexity of hafnians is understood less well than that of permanents. Unlike in the case of the permanent, no fully polynomial (randomized or deterministic) polynomial approximation scheme is known to compute the hafnian of a non-negative real symmetric matrix. Unlike in the case of the permanent, no deterministic polynomial time algorithm approximating the hafnian of a $2n \times 2n$ non-negative symmetric matrix within a factor of $c^n$, where $c>0$ is an absolute constant, is known.
On the other hand there is a polynomial time randomized algorithm based on the representation of the hafnian as the expectation of the determinant of a random matrix, which approximates the hafnian of a given non-negative symmetric $2n \times 2n$ matrix within a factor of $c^n$, where $c \approx 0.56$ \cite{Ba99}. Also, for any $0 < \epsilon < 1$ fixed in advance, there is a deterministic polynomial time algorithm based on scaling, which, given a $2n \times 2n$ symmetric matrix 
$A=\left(a_{ij}\right)$ satisfying 
$$\epsilon \ \leq \ a_{ij} \ \leq \ 1 \quad \text{for all} \quad i,j,$$
computes $\haf A$ within a multiplicative factor of $n^{\kappa(\epsilon)}$ for some $\kappa(\epsilon)>0$ \cite{BS11}.

With minimal changes, the approach of this paper can be applied to computing hafnians. Namely, let $J$ denote the $2n \times 2n$ matrix filled with 1s and let us define 
$$f(z)=\ln \haf \bigl(J + z(A-J)\bigr).$$
Then 
$$f(0)=\ln \haf J = \ln {(2n)! \over n! 2^n} \quad \text{and} \quad f(1)= \ln \haf A$$ and one can use the Taylor polynomial approximation (1.1.2) to estimate $f(1)$. As in Section 2, one can compute the right hand side of (1.1.2) in $n^{O(m)}$ time.
The statement and the proof of Theorem 1.3 carries over to hafnians almost verbatim. Namely, let $\delta>0$ be a real for which the equation (3.4.1) has a solution $0 < \theta < \pi/2$ (hence one can choose $\delta=0.195$). Then 
$\haf Z \ne 0$ as long as $Z=\left(z_{ij}\right)$ is a $2n \times 2n$ symmetric complex matrix satisfying 
$$\left| z_{ij}-1\right| \ \leq \ \delta \quad \text{for all} \quad i,j.$$
Instead of the row expansion of the permanent (3.2.1) used in Lemmas 3.2 and 3.3, one should use the row expansion of the hafnian
$$\haf Z=\sum_{j=2}^{2n} z_{1j} \haf Z_j,$$ 
where $Z_j$ is the symmetric $(2n-2) \times (2n-2)$ matrix obtained from $Z$ by crossing out the first and the $j$-th row and the first and the $j$-th column. As in Section 2, we obtain an algorithm of $n^{O\left(\ln n - \ln \epsilon\right)}$ complexity of approximating $\haf Z$ within relative error $\epsilon>0$, where $Z=\left(Z_{ij}\right)$ is a $2n \times 2n$ symmetric complex matrix satisfying 
$$\left| z_{ij}-1\right| \ \leq \ \gamma, \quad \text{for all} \quad i, j.$$
and $\gamma > 0$ is an absolute constant (one can choose $\gamma=0.19$).

\subhead (4.2) Multidimensional permanents \endsubhead Let us fix an integer $\nu \geq 2$ and let 
$$A=\left(a_{i_1 \ldots i_{\nu}}\right), \quad 1 \leq i_1, \ldots, i_{\nu} \leq n,$$ be an 
$\nu$-dimensional cubical $n \times \ldots \times n$ array of real or complex numbers. We define 
$$\PER A =\sum_{\sigma_1, \ldots, \sigma_{\nu-1} \in S_n} \prod_{i=1}^n a_{i \sigma_1(i) \ldots \sigma_{\nu-1}(i)}.$$
If $\nu=2$ then $A$ is an $n \times n$ matrix and $\PER A=\per A$. For $\nu >2$ it is already an NP-hard problem to tell 
$\PER A$ from $0$ even if $a_{i_1 \ldots i_{\nu}} \in \{0, 1\}$ since the problem reduces to detecting a perfect matching in a hypergraph, see, for example, Problem SP1 in \cite{A+99}. However, for any $0 < \epsilon < 1$, fixed in advance, there is a polynomial time deterministic algorithm based on scaling, which, given a real array $A$ satisfying
$$\epsilon \ \leq \ a_{i_1 \ldots i_{\nu}} \ \leq \ 1 \quad \text{for all} \quad 1 \leq i_1, \ldots, i_{\nu} \leq n$$
computes $\PER A$ within a multiplicative factor of $n^{\kappa(\epsilon, \nu)}$ for some $\kappa(\epsilon, \nu) >0$ \cite{BS11}.

With some modifications, the method of this paper can be applied to computing this multidimensional version of the permanent.
Namely, let $J$ be the array filled with 1s and let us define 
$$f(z)=\ln \PER \bigl(J +z(A-J)\bigr).$$
Then 
$$f(0)=\ln \PER J = (\nu-1) \ln n! \quad \text{and} \quad f(1)=\ln \PER A$$
and one can use the Taylor polynomial approximation (1.1.2) to estimate $f(1)$. As in Section 2, one can compute the right hand side of (1.1.2) in $n^{O(m)}$ time, where the implicit constant in ``$O(m)$" depends on $\nu$. The proof of Theorem 1.3 carries to multidimensional permanents with some modifications. Namely, for some sufficiently small $\delta_{\nu}>0$ the equation 
$$\theta={2 \delta_{\nu} \over \left(1-\delta_{\nu}\right) \sqrt{\cos\bigl((\nu-1) \theta + 2 \arcsin \delta_{\nu}\bigr)}}$$
has a solution $\theta \geq 0$ such that $(\nu-1) \theta + 2 \arcsin \delta_{\nu} < \pi/2$. For $\nu=2$, we get the equation (3.4.1) with a possible choice of $\delta_2=0.195$, while for 
$\nu=3$ we can choose $\delta_3=0.125$ and for $\nu=4$ we can choose $\delta_4=0.093$. Then $\PER Z \ne 0$ as long as 
$Z=\left(z_{i_1 \ldots i_{\nu}}\right)$ is an array of complex numbers satisfying 
$$\left| z_{i_1 \ldots i_{\nu}}-1\right| \ \leq \ \delta_{\nu} \quad \text{for all} \quad 1 \leq i_1, \ldots, i_{\nu} \leq n.$$
We proceed as in the proof of Theorem 1.3, only instead of the first row expansion of the permanent (3.2.1) used in Lemmas 3.2 and 3.3, we use the first index expansion
$$\PER Z =\sum_{1 \leq j_2, \ldots, j_{\nu} \leq n} z_{1 j_2 \ldots j_{\nu}} \PER Z_{j_2 \ldots j_{\nu}},$$
where $Z_{j_2 \ldots j_{\nu}}$ is the $\nu$-dimensional array of size $(n-1) \times \cdots \times (n-1)$ obtained from $Z$ by crossing out the section with the first index 1, the section with the second index $j_2$ and so forth, concluding with crossing out the section with the last index $j_{\nu}$. As in Section 2, we obtain at algorithm of $n^{O(\ln n -\ln \epsilon)}$ complexity of approximating $\PER Z$ within relative error $\epsilon >0$, where $Z$ is a $\nu$-dimensional cubic $n \times \cdots \times n$ array of complex numbers satisfying
$$\left| z_{i_1 \ldots i_{\nu}} -1\right| \ \leq \ \gamma_{\nu}  \quad \text{for all} \quad 1\leq i_1, \ldots, i_{\nu} \leq n,$$
and $0< \gamma_{\nu} < \delta_{\nu}$ are absolute constants (one can choose $\gamma_2=0.19$, $\gamma_3=0.12$ and 
$\gamma_4=0.09$).

\enddocument

\end